# Self-optimal Piecewise Linearization Based Network Power Flow Constraints in Electrical Distribution System Optimization

Yun Zhou, *Member, IEEE*, Jingyang Yun, *Student Member, IEEE*, Zheng Yan, Donghan Feng, *Senior Member, IEEE*, and Sijie Chen, *Member, IEEE*

*Abstract*—As a representative mathematical expression of power flow (PF) constraints in electrical distribution system (EDS), the piecewise linearization (PWL) based PF constraints have been widely used in different EDS optimization scenarios. However, the linearized approximation errors originating from the currently-used PWL approximation function can be very large and thus affect the applicability of the PWL based PF constraints. This letter analyzes the approximation error self-optimal (ESO) condition of the PWL approximation function, refines the PWL function formulas, and proposes the self-optimal (SO)-PWL based PF constraints in EDS optimization which can ensure the minimum approximation errors. Numerical results demonstrate the effectiveness of the proposed method.

## I. Introduction

The Distflow network PF constraints linearized with PWL approximation function constitute the PWL based PF constraints [1], and have been widely used in different EDS optimization scenarios. This letter analyzes the causes of linearized approximation errors in the PWL approximation function, and proposes the ESO condition of the PWL function. The currently-used PWL function formulas are complemented to meet the ESO condition. And the SO-PWL based network PF constraints are realized by implementing the supplemented PWL function in the initial PWL based PF constraints. The service restoration problem in EDS is selected to demonstrate the effectiveness of the proposed method in numerical results.

## II. Fundamental PWL Approximation Procedure of the Distflow Network PF Constraints

In the Distflow network PF constraints [1], only the constraint shown in (1) contains quadratic variables, and is nonlinear.

$$(V_{\text{bus}}^{\text{norm}})^2 I_{ij}^{\text{sqr}} = P_{ij}^2 + Q_{ij}^2 \quad (1)$$

$$P_{ij}^2 \approx f(P_{ij}, P_{ij}^{\max}, \Lambda) \quad (2)$$

$$Q_{ij}^2 \approx f(Q_{ij}, Q_{ij}^{\max}, \Lambda) \quad (3)$$

$$(V_{\text{bus}}^{\text{norm}})^2 I_{ij}^{\text{sqr}} = f(P_{ij}, P_{ij}^{\max}, \Lambda) + f(Q_{ij}, Q_{ij}^{\max}, \Lambda) \quad (4)$$

The PWL approximation procedure of the Distflow network PF constraints concentrates mainly on the linearization of constraint (1) (shown in (2)-(4)). The currently-used fundamental PWL approximation function formulas are shown in (5)-(12) [1].

$$f(y, \bar{y}, \Lambda) = \sum_{\lambda=1}^{\Lambda} \phi_{y,\lambda} \Delta_{y,\lambda} \quad (5)$$

$$y = y^+ - y^- \quad (6)$$

$$y^+ + y^- = \sum_{\lambda=1}^{\Lambda} \Delta_{y,\lambda} \quad (7)$$

$$0 \leq \Delta_{y,\lambda} \leq \bar{y}/\Lambda \quad \forall \lambda=1,2,...,\Lambda \quad (8)$$

$$\phi_{y,\lambda} = (2\lambda-1)\bar{y}/\Lambda \quad \forall \lambda=1,2,...,\Lambda \quad (9)$$

$$0 \leq y^+ \leq z^+ \bar{y} \quad (10)$$

$$0 \leq y^- \leq z^- \bar{y} \quad (11)$$

$$z^+ + z^- \leq 1 \quad (12)$$

In (12), $z^+$ and $z^-$ are binary variables. The definitions of variables in (1)-(12) can be referred to [1]. Constraint (4) and other linear constraints in the Distflow network PF constraints compose the initial PWL based network PF constraints [1].

## III. ESO Condition of the PWL Approximation Function and SO-PWL Based Network PF Constraints

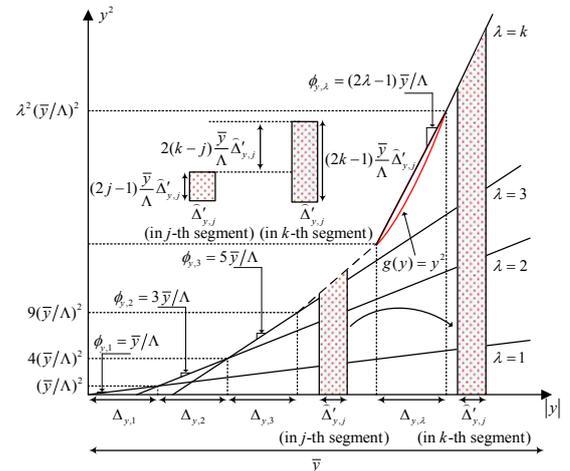

Fig. 1. Schematic diagram of the PWL function approximating the quadratic curve

For $g(y) = y^2 \approx f(y, \bar{y}, \Lambda)$, suppose $y \geq 0$ for convenience, the relative PWL approximation error $E_y$ is defined in (13).

$$E_y = |f(y, \bar{y}, \Lambda) - y^2| / y^2 \times 100\% \quad (13)$$

$$\lambda_{y,\text{up}} = \begin{cases} 1 & y = 0 \\ \text{ceil}(y/(\bar{y}/\Lambda)) & y > 0 \end{cases} \quad (14)$$

$$C_{\text{ESO}} = \left\{ \Delta_{y,\lambda} \middle| \begin{array}{l} \Delta_{y,\lambda} = \overline{y}/\Lambda \quad \lambda = 1, 2, ..., \lambda_{y,\text{up}} - 1, \lambda_{y,\text{up}} > 1 \\ 0 \le \Delta_{y,\lambda} \le \overline{y}/\Lambda \quad \lambda = \lambda_{y,\text{up}} \\ \Delta_{y,\lambda} = 0 \quad \lambda = \lambda_{y,\text{up}} + 1, \lambda_{y,\text{up}} + 2, ..., \Lambda \end{array} \right\} \quad (15)$$

To minimize the PWL approximation error, the ESO condition of the PWL approximation function is defined in (15), in which $\lambda_{y,\text{up}}$ is defined in (14). The ESO condition in (15) implies that the segments $\Delta_{y,\lambda}$ must be used (i.e. take a certain value) one after another during approximation process. And the segments in use state except the last one should equal to its upper limit $\overline{y}/\Lambda$. The following propositions and proofs illustrate the validity of the ESO condition:

*1) Proposition 1*: The PWL approximation results is no less than the real variable values (i.e. $f(y,\overline{y},\Lambda) \ge y^2$).

*Proof 1*: Fig. 1 shows how the PWL function approximates the quadratic curve $g(y)$. Referring to (5)-(12), due to the convexity of $g(y)$, all segments employed in the PWL function are above on $g(y)$ (shown in Fig. 1), and the Proposition 1 is correct obviously.

*2) Proposition 2*: The ESO condition is sufficient to achieve the minimum PWL approximation error.

*Proof 2*: Let $\{\Delta_{y,\lambda}\}$ be the set of segments derived from (5)-(12), and satisfy the ESO condition in (15). Suppose $\{\Delta'_{y,\lambda}\}$ be another set of segments not satisfying (15). Referring to (15), for $\lambda_{y,\text{up}} > 1$, there must be at least one segment in the preceding $\lambda_{y,\text{up}} - 1$ segments in $\{\Delta'_{y,\lambda}\}$ failing to fill the full interval (e.g. the *j*-th segment $\Delta'_{y,j}$ in (16)). And the remainder term of $\Delta'_{y,j}$ (i.e. $\hat{\Delta}'_{y,j}$ in (17)) should be compensated in certain later segments in $\{\Delta'_{y,\lambda}\}$. (e.g. compensated in the *k*-th segment in $\{\Delta'_{y,\lambda}\}$, shown in (18) and Fig. 2).

$$0 \le \Delta'_{y,j} < \overline{y}/\Lambda \quad j \in \{1, 2, ..., \lambda_{y,\text{up}} - 1\} \quad (16)$$

$$\hat{\Delta}'_{y,j} = \overline{y}/\Lambda - \Delta'_{y,j} \quad (17)$$

$$\Delta'_{y,k} - \Delta_{y,k} \ge \hat{\Delta}'_{y,j} \quad k \in \{\lambda_{y,\text{up}}, \lambda_{y,\text{up}} + 1, ..., \Lambda\} \quad (18)$$

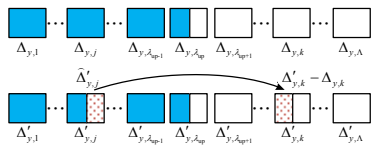

Fig. 2. Schematic diagram of filling states of $\{\Delta_{y,\lambda}\}$ and $\{\Delta'_{y,\lambda}\}$

$$\sum_{\lambda=1}^{\Lambda} \phi_{y,\lambda} \Delta'_{y,\lambda} - y^2 \ge \sum_{\lambda=1}^{\Lambda} \phi_{y,\lambda} \Delta_{y,\lambda} - y^2 \ge 0 \quad (19)$$

The increment in $f(y,\overline{y},\Lambda)$ due to the compensation of $\hat{\Delta}'_{y,j}$ is shown in Fig. 1. As the slope of segment (i.e. $\phi_{y,\lambda}$) is monotonic increasing with $\lambda$, with Proposition 1 and $E_y$ defined in (13), (19) is correct and Proposition 2 is verified.

For some optimization scenarios with particular objective (e.g. containing the minimum power loss item in the objective function), the ESO condition can be ensured implicitly [2]. For more general scenarios, the PWL function formulas need to be complemented to realize the ESO condition.

$$\Delta_{y,\lambda} - \overline{y}/\Lambda + (1 - x_\lambda)M + \varepsilon^+ \ge 0 \quad \forall \lambda = 1, 2, ..., \Lambda \quad (20)$$

$$0 \le \Delta_{y,\lambda+1} \le x_\lambda \overline{y}/\Lambda \quad \forall \lambda = 1, 2, ..., \Lambda - 1 \quad (21)$$

By introducing a sufficiently big positive constant $M$, a near-zero positive constant $\varepsilon^+$, and a set of binary variables $\{x_\lambda\}$, the fundamental PWL approximation function can be supplemented as (5)-(12), (20), and (21) to meet the ESO condition. And by implementing the supplemented PWL function in the initial PWL based network PF constraints, the SO-PWL based network PF constraints are realized.

## IV. NUMERICAL RESULTS

The EDS service restoration is to restore the maximum outage loads with limited generating capacity (i.e. objective function in (22) with PF constraints as major constraints).

$$\max \sum_{i \in S_{\text{bus}}} P_{\text{L},i} \quad (22)$$

The proposed method is tested on the typical IEEE 33-bus distribution system plus 4 DGs (located at bus 13, 21, 22 and 30 with maximum active/reactive power limit of 0.05/0.03 pu each), $I_{ij}^{\max} = 50$ A and $\Lambda = 50$. Table I tabulates the comparison of $E_p$ of $P_{ij}^2$ in 3 feeders (feeders of relative bigger $E_p$ in $S_{\text{feeder}}$) in the service restoration solving results. The corresponding filling states of $\{\Delta_{p,\lambda}\}$ are shown in Fig.3.

TABLE I
THE PWL APPROXIMATION ERRORS OF $P_{ij}^2$ ($E_p$) OF FEEDERS WITH PWL/SO-PWL BASED PF CONSTRAINTS

| Feeder | $E_p$ (%) | |
|---|---|---|
| | PWL | SO-PWL |
| 9-10 | 241.310 | 0.848 |
| 12-13 | 68.407 | 0.469 |
| 13-14 | 62.520 | 0.519 |

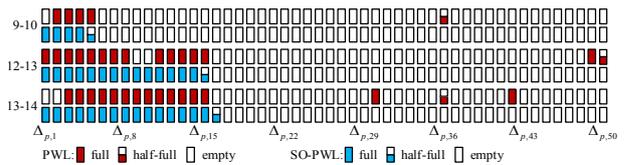

Fig. 3. Schematic diagram of filling states of $\{\Delta_{p,\lambda}\}$ in the solving results

From Fig. 3 and Table I, for the PWL based network PF constraints case, $\{\Delta_{p,\lambda}\}$ are filled unorderly with large $E_p$. For the SO-PWL based PF constraints case, the filling states are well-ordered (satisfying the ESO condition) and the $E_p$ meet significant reduction, therefore the effectiveness of the proposed SO-PWL based PF constraints can be demonstrated.